\documentclass{article}
\usepackage[utf8]{inputenc}
\usepackage{amsmath,amsthm,amssymb}
\usepackage{mathrsfs,comment}
\usepackage{mathtools}
\usepackage{xspace}

\newcommand{\C}{\mathbb{C}}

\newcommand{\SO}{\text{SO}}

\newtheorem{theorem}{Theorem}[section]

\theoremstyle{definition}

\theoremstyle{remark}

\title{Diagonalization of the metric of a Lorentzian 3-manifold}
\author{Romeo Segnan Dalmasso}

\begin{document}

\maketitle

\begin{abstract}
We study the problem of diagonalization of the metric of 3-dimensional Lorentzian manifold. Applying the technique of moving frames, we prove that every smooth Lorentzian 3-manifold admits an atlas in which the metric assumes a diagonal form.
\end{abstract}

\renewcommand{\thefootnote}{\fnsymbol{footnote}}
\footnotetext{\emph{MSC class 2020}: \emph{Primary} 53B30; \emph{Secondary} 58J60
}
\renewcommand{\thefootnote}{\arabic{footnote}}

\section{Introduction}
A (pseudo)-Riemannian $n$-manifold $(M,g)$ is said to have \emph{orthogonal coordinates around a point} if in the neighborhood of the point there is a chart such that the metric $g$ with respect to it is in diagonal form, i.e. 
\[
g=\sum_{i=1}^nf_idx^i\otimes dx^i.
\]
If the manifold satisfies this property at every point, then we will say that it admits an \emph{orthogonal atlas}. 

Any surface has an orthogonal atlas since there always are isothermal coordinates (see \cite{Deturck1981SomeGeometry}) and the metric assumes the particular diagonal form 
\[
g=f(x,y)\big(dx\otimes dx+dy\otimes dy\big).
\]
In the Riemannian setting D. DeTurck and D. Yang in \cite{Deturck1984} proved that any 3-dimensional smooth manifold $(M,g)$ has an orthogonal atlas using the technique of moving frames. In this case the metric assumes the more general form
\begin{equation*}
    g=f_1(x, y,z)dx\otimes dx+f_2(x, y,z)dy\otimes dy+f_3(x, y,z)dz\otimes dz.
\end{equation*}
In their paper they point out that in higher dimension the situation changes because the existence of the orthogonal atlas is subject to a condition on the Weyl tensor. Subsequently, P. Tod in \cite{Tod1992OnCC} studied the problem in the same setting but in dimension $n\ge4$ where he found necessary algebraic conditions for the existence of the orthogonal atlas. In the paper the cases of dimensions $n=4$, $n=5$ and $n\ge6$ are studied separately as they require each a different condition on the Weyl tensor of the manifold; more precisely  the restriction involves the third derivatives of the tensor for $n=4$, the first derivative for $n=5$,  and the tensor itself for $n\ge6$. However, J. Grant and J. A. Vickers in \cite{Grant2009BlockMetrics} proved that something can be said even in dimension 4, in particular they showed that in the analytic setting, in the Riemannian and in Lorentzian case, one can find a chart such that the metric $g$ is block diagonal, i.e.
\[
g_{ij}=\begin{pmatrix} A & 0 \\ 0 & B \end{pmatrix}
\]
where $A$ and $B$ are 2×2 block matrix, even when no assumptions are made on the Weyl tensor (or its derivatives).

More recently, O. Kowalski and M. Sekizawa in \cite{KOWALSKI2013251} proved the existence of an orthogonal atlas in the real analytic Lorentzian setting by applying the Cauchy-Kovalevski Theorem. On the other hand, P. Gauduchon and A. Moroianu proved in \cite{Gauduchon2020Non-existenceSpaces} that one cannot find orthogonal coordinates in the neighborhood of any point for the complex and quaternionic projective spaces $\C\mathbb{P}^m$ and $\mathbb{HP}^q$. 

In this paper we will prove that all smooth Lorentzian 3-manifolds admit an orthogonal atlas by following the same method used by DeTurck and Yang. The technique is the following: the problem is initially shifted from a PDE system about the coordinates to a PDE system about a coframe in the cotangent bundle with respect to a fixed coframe, then one proves that the Cauchy problem associated to the second PDE system admits a solution. 

\textbf{Acknowledgments}
This paper was written as part of the author's PhD thesis for the joint PhD program in Mathematics Università di Milano Bicocca - University of Surrey, and it is based on a suggestion of his supervisor James Grant. The author also thanks his other supervisor Diego Conti for some useful advice.

The author acknowledges GNSAGA of INdAM.

\section{Orthogonal coordinates on Lorentzian manifolds}
We proceed to illustrate the proof of the following
\begin{theorem}\label{teo:main}
Let $(M,g)$ be a smooth Lorentzian 3-manifold. Then $M$ admits an orthogonal atlas.
\end{theorem}
Let $(\bar e_1,\bar e_2,\bar e_3)$ be an orthonormal frame of $(M,g)$ and $(\bar\omega_1,\bar\omega_2,\bar\omega_3)$ the corresponding coframe. We want to find a triplet of coordinated functions $(x_1,x_2,x_3)$ such that, if $e_i=\partial_i$ is the coordinated frame of $(x_1,x_2,x_3)$, then $g(e_i,e_j)=0$ every time $i\ne j$. The first difficulty we find both in the Riemannian and in the Lorentzian setting is the following. Assume $(y^1,y^2,y^3)$ are fixed coordinates and $g(y)$ is the metric tensor w.r.t. this chart; then the coordinated frame $\{e_i\}$ can be written as
\[
e_i=\frac{\partial y^\alpha}{\partial x^i}\frac{\partial}{\partial y^\alpha}
\]
and hence the PDE system to be solved is
\begin{equation}
    0=g(\partial_i,\partial_j)=\sum_{\alpha,\beta=1}^3\frac{\partial y^\alpha}{\partial x^i}\,\frac{\partial y^\beta}{\partial x^j}\,g_{\alpha\beta}(y)\quad\text{for }i\ne j.
\end{equation}
This system is nonlinear, and its linearization is not symmetric hyperbolic, which means that the standard results of existence of the solution do not apply. Furthermore, there is an invariance in the solution if the unknowns are the coordinates: assume $(\tilde x^1,\tilde x^2,\tilde x^3)$ are other coordinates, such that $\tilde x^i=f^i(x^i)$ and each $f^i$ is a strictly monotone function. Then
\[
0=g\left(\frac{\partial}{\partial x^i},\frac{\partial}{\partial x^j}\right)=\frac{\partial f^i}{\partial x^i}\frac{\partial f^j}{\partial x^j}g\left(\frac{\partial}{\partial \tilde x^i},\frac{\partial}{\partial \tilde x^j}\right)
\]
and hence also $(\tilde x^1,\tilde x^2,\tilde x^3)$ are orthogonal coordinates.

For this reason it works best if one does not set the unknowns to be the coordinated functions $(x_1,x_2,x_3)$, but the normalized coframe $(\omega^1,\omega^2,\omega^3)$, where $\omega^i=f_idx^i$ (no sum intended) and $f_i=1/|dx^i|$. Applying Frobenius Theorem it is easy to get an equivalent condition to the existence of the coordinated charts depending on the coframe, that is
\begin{equation}\label{eqn:FrobeniusCondition}
\omega^i\wedge d\omega^i=0
\end{equation}
must hold, when $i=1,2,3$. Now, as the coframe has to be orthonormal, it has to satisfy the first Cartan structure equation
\[
d\omega^i=\sum_j\omega^j\wedge\omega_j^i
\]
where $(\omega^j_i)$ is the connection matrix 1-form. Here appears the first difference between the Riemannian and Lorentzian case, although it does not yield any actual change in the proof: in the first case $\omega_i^j=-\omega_j^i$ for any $i,j$, but in the second we have 
\begin{equation}\label{eqn:firstDifference}
    \omega_1^2=-\omega_2^1,\quad\omega_1^3=\omega_3^1,\quad\omega_2^3=\omega_3^2.
\end{equation}
Hence \eqref{eqn:FrobeniusCondition} becomes
\begin{align*}
    \omega^1\wedge\omega^2\wedge\omega_2^1+\omega^1\wedge\omega^3\wedge\omega_3^1&=0\\
    \omega^1\wedge\omega^2\wedge\omega_2^1+\omega^2\wedge\omega^3\wedge\omega_3^2&=0\\
    \omega^1\wedge\omega^3\wedge\omega_3^1+\omega^2\wedge\omega^3\wedge\omega_3^2&=0,
\end{align*}
thus, by alternatively subtracting one and adding the other we get the system
\begin{equation}\label{eqn:equationSystem}
    \omega^1\wedge\omega^2\wedge\omega_2^1=0,\quad\omega^2\wedge\omega^3\wedge\omega_3^2=0,\quad\omega^1\wedge\omega^3\wedge\omega_3^1=0.
\end{equation}
We now write $\omega^i$ with respect to $\bar\omega^j$ and vice-versa as
\[
\omega^i=b^i_j\bar\omega^j,\quad\bar\omega^j=\bar b^i_j\omega^i
\]
and we will solve for the $b^j_i$. We will solve \eqref{eqn:equationSystem}, hence we need $\omega^i_j$ and we start by noting that
\[
\begin{split}
    \omega^l\wedge\omega^i_l=d\omega^i&=d\sum_jb^i_j\bar\omega^j=\sum_j\left(\sum_k\bar e_k(b^i_j)\bar\omega^k\wedge\bar\omega^j+b^i_j\bar\omega^k\wedge\bar\omega^i_k\right)\\
    &=\sum_{j,k}\bar\omega^k\wedge\big(\bar e_k(b^i_j)\bar\omega^j+b^i_j\bar\omega^i_k\big)\\
    &=\sum_{j,k,l}\omega^l\wedge\big(b^l_k\bar e_k(b^i_j)\bar\omega^j+b^l_kb^i_j\bar\omega^i_k\big).
\end{split}
\]
As a consequence of the first difference we find a second one here: while
\[
\omega^1_2=\sum_{j,k}\frac{1}{2}\big\{b^2_k\bar e_k(b^1_j)-b^1_k\bar e_k(b^2_j)\big\}\bar\omega^j+b^2_kb^i_j\bar\omega^1_k
\]
remains as in \cite{Deturck1984}, the other two differ due to \eqref{eqn:firstDifference} as follows:
\[
\omega^1_3=\sum_{j,k}\frac{1}{2}\big\{b^3_k\bar e_k(b^1_j)+b^1_k\bar e_k(b^3_j)\big\}\bar\omega^j+b^3_kb^i_j\bar\omega^1_k
\]
and
\[
\omega^2_3=\sum_{j,k}\frac{1}{2}\big\{b^3_k\bar e_k(b^2_j)+b^2_k\bar e_k(b^3_j)\big\}\bar\omega^j+b^3_kb^i_j\bar\omega^1_k.
\]
Again by following \cite{Deturck1984} we rewrite \eqref{eqn:equationSystem} substituting $\omega^i_j$ and obtaining
\begin{gather*}
0=\sum_{i,l,j,k}b^1_ib^2_l\bar\omega^i\wedge\bar\omega^l\wedge\left[\frac{1}{2}\big\{b^2_k\bar e_k(b^1_j)-b^1_k\bar e_k(b^2_j)\big\}\bar\omega^j+b^2_kb^i_j\bar\omega^1_k\right]\\
0=\sum_{i,l,j,k}b^1_ib^3_l\bar\omega^i\wedge\bar\omega^l\wedge\left[\frac{1}{2}\big\{b^3_k\bar e_k(b^1_j)+b^1_k\bar e_k(b^3_j)\big\}\bar\omega^j+b^3_kb^i_j\bar\omega^1_k\right]\\
0=\sum_{i,l,j,k}b^2_ib^3_l\bar\omega^i\wedge\bar\omega^l\wedge\left[\frac{1}{2}\big\{b^3_k\bar e_k(b^2_j)+b^2_k\bar e_k(b^3_j)\big\}\bar\omega^j+b^3_kb^i_j\bar\omega^1_k\right]
\end{gather*}
The unknowns of the system are $(b_i^j)\in C^\infty(M,\SO(2,1))$.

We are going to prove that the linearization of this system is diagonal hyperbolic. Consider the linearization $\beta_j^i=(\delta b)_j^i$ and notice that we can assume that $\{\bar\omega^i\}=\{\omega^i\}$ when we linearize around $\{\omega^i\}$, as such $b^i_j(x)=\delta^i_j$. Thus, the linearized system is
\begin{gather*}
0=\delta^1_i\delta^2_l\frac{1}{2}\big(\delta^2_k\bar e_k(\beta^1_j)-\delta^1_k\bar e_k(\beta^2_j)\big)\bar\omega^i\wedge\bar\omega^l\wedge\bar\omega^j+\text{ lower order terms in }\beta\\
0=\delta^1_i\delta^3_l\frac{1}{2}\big(\delta^3_k\bar e_k(\beta^1_j)+\delta^1_k\bar e_k(\beta^3_j)\big)\bar\omega^i\wedge\bar\omega^l\wedge\bar\omega^j+\text{ lower order terms in }\beta\\
0=\delta^2_i\delta^3_l\frac{1}{2}\big(\delta^3_k\bar e_k(\beta^2_j)+\delta^2_k\bar e_k(\beta^3_j)\big)\bar\omega^i\wedge\bar\omega^l\wedge\bar\omega^j+\text{ lower order terms in }\beta
\end{gather*}
in which the only non-zero elements are
\begin{gather*}
\frac{1}{2}\big(\bar e_2(\beta^1_3)-\bar e_1(\beta^2_3)\big)=\text{ terms of order 0 in }\beta,\\
\frac{1}{2}\big(\bar e_3(\beta^1_2)+\bar e_1(\beta^3_2)\big)=\text{ terms of order 0 in }\beta,\\
\frac{1}{2}\big(\bar e_3(\beta^2_1)+\bar e_2(\beta^3_1)\big)=\text{ terms of order 0 in }\beta.
\end{gather*}
As $(b^i_j(x))\in\SO(2,1)$ we have that $(\beta^i_j)\in\mathfrak{so}(2,1)$, hence we can rewrite everything as
\begin{gather*}
\bar e_1(\beta^2_3)=\text{ terms of order 0 in }\beta\\
\bar e_2(\beta^1_3)=\text{ terms of order 0 in }\beta\\
\bar e_3(\beta^1_2)=\text{ terms of order 0 in }\beta.
\end{gather*}
The differential operator is thus 
\[
A(u)=\bar e_1(u)+\bar e_2(u)+\bar e_3(u)
\]
that is in diagonal form, and its symbol is
\[
\sigma(\xi)=\sum_{i=1}^3\xi^i.
\]
To finally prove that the metric is diagonalizable we have to find a solution to the Cauchy problem given by the differential operator $A$ and a set of initial data to be chosen. To do so, we need these data to not be characteristic of the operator. By the form of $A$ we deduce that the characteristics of the system are the covectors that annihilate $e_1,e_2$ and $e_3$. Hence, the initial data for the Cauchy problem associated to the system can be given as  the coframe $\{\omega^i\}$ on a surface $\Sigma\subset M$ with $e_i\notin T\Sigma$ for $i=1,2,3$ since we need $\omega^i(v)\ne0$ for all $v\in T\Sigma$. This concludes the proof of Theorem \ref{teo:main}.

\bibliographystyle{acm}
\bibliography{references}

\begin{thebibliography}{1}

\bibitem{Deturck1981SomeGeometry}
{\sc Deturck, D.~M., and Kazdan, J.~L.}
\newblock {Some regularity theorems in Riemannian geometry}.
\newblock {\em Annales scientifiques de l'{\'{E}}cole normale sup{\'{e}}rieure
  14}, 3 (1981), 249--260.

\bibitem{Deturck1984}
{\sc DeTurck, D.~M., and Yang, D.}
\newblock {Existence of elastic deformations with prescribed principal strains
  and triply orthogonal systems}.
\newblock {\em Duke Mathematical Journal 51}, 2 (1984), 243--260.

\bibitem{Gauduchon2020Non-existenceSpaces}
{\sc Gauduchon, P., and Moroianu, A.}
\newblock {Non-existence of orthogonal coordinates on the complex and
  quaternionic projective spaces}.
\newblock {\em Journal of Geometry and Physics 155\/} (9 2020), 103770.

\bibitem{Grant2009BlockMetrics}
{\sc Grant, J. D.~E., and Vickers, J.~A.}
\newblock {Block diagonalization of four-dimensional metrics}.
\newblock {\em Classical and Quantum Gravity 26}, 23 (12 2009), 235014.

\bibitem{KOWALSKI2013251}
{\sc Kowalski, O., and Sekizawa, M.}
\newblock {Diagonalization of three-dimensional pseudo-Riemannian metrics}.
\newblock {\em Journal of Geometry and Physics 74\/} (2013), 251--255.

\bibitem{Tod1992OnCC}
{\sc Tod, K.~P.}
\newblock {On choosing coordinates to diagonalize the metric}.
\newblock {\em Classical and Quantum Gravity 9\/} (1992), 1693--1705.

\end{thebibliography}

\end{document}